\documentclass[12pt,a4paper]{article}
\usepackage[width=18cm,height=25cm]{geometry}
\usepackage{amsmath,amssymb,amsfonts,amsthm,enumerate,float}

\usepackage{bbm}
\usepackage{graphicx}
\usepackage{tikz}
\usetikzlibrary{shapes,backgrounds}
\usepackage{tikz-cd}
\usetikzlibrary{%
  matrix,%
  calc,%
  arrows%
}

\usepackage{scalerel}

\usepackage{mathrsfs,amsbsy,bm,mathtools}
\usepackage{indent first}

\newcommand{\bc}{\mathbb{C}}

\newcommand{\br}{\mathbb{R}}

\newtheorem{theorem}{Theorem}[section]
\newtheorem{lemma}[theorem]{Lemma}

\newtheorem{Prop}[theorem]{Proposition}
\newtheorem{Cor}[theorem]{Corollary}
\newtheorem{Quest}[theorem]{Question}

\theoremstyle{definition}
\newtheorem{definition}[theorem]{Definition}

\title{\bf Hodge-Riemann property of Griffiths positive matrices with $(1,1)$-form entries}
\author{Zhangchi Chen}
\date{\today}

\begin{document}

\maketitle

\begin{abstract}
The classical Hard Lefschetz theorem (HLT), Hodge-Riemann bilinear relation theorem (HRR) and Lefschetz decomposition theorem (LD) are stated for a power of a K\"ahler class on a compact K\"ahler manifold. These theorems are not true for an arbitrary class, even if it contains a smooth strictly positive representative.

Dinh-Nguy\^en proved the mixed HLT, HRR and LD for a product of arbitrary K\"ahler classes. Instead of products, they asked whether determinants of Griffiths positive $k\times k$ matrices with $(1,1)$-form entries in $\bc^n$ satisfies these theorems in the linear case.

This paper answered their question positively when $k=2$ and $n=2,3$. Moreover, assume that the matrix only has diagonalized entries, for $k=2$ and $n\geqslant 4$, the determinant satisfies HLT for bidegrees $(n-2,0)$, $(n-3,1)$, $(1,n-3)$ and $(0,n-2)$. In particular, for $k=2$ and $n=4,5$ with this extra assumption, the determinant satisfies HRR, HLT and LD.

Two applications: First, a Griffiths positive $2\times 2$ matrix with $(1,1)$-form entries, if all entries are $\mathbb{C}$-linear combinations of the diagonal entries, then its determinant also satisfies these theorems. Second, on a complex torus of dimension $\leqslant 5$, the determinant of a Griffiths positive $2\times 2$ matrix with diagonalized entries satisfies these theorems.

\end{abstract}

\section{Introduction}
\label{sect-intro}
Let $X$ be a compact K\"ahler manifold of dimension $n$ and let $\omega$ be a K\"ahler form. The cohomology of $X$ satisfies the Hodge-decomposition:
\[
H^d(X,\bc)=\bigoplus\limits_{p+q=d}H^{p,q}(X,\bc), \ \ \ (0\leqslant d\leqslant n); \ \ \ \ \overline{H^{p,q}(X,\bc)}=H^{q,p}(X,\bc);
\]
where $H^{p,q}(X,\bc)$ is the Hodge cohomology group of bidegree $(p,q)$ of $X$ with the convention that $H^{p,q}(X,\bc)=0$ unless $0\leqslant p,q\leqslant n$. When $p,q\geqslant0$ and $p+q\leqslant n$, let $k=n-p-q$ and define $\Omega:=\omega^k$ a $(k,k)$-form on $X$.

Recall the classical hard Lefschetz theorem (HLT), the Hodge-Riemann bilinear relation theorem (HRR) and the Lefschetz decomposition theorem (LD). One may refer to BDIP \cite{BDIP-1996}, Griffiths and Harris \cite{Griffiths-Harris-1994} and Voisin \cite{Voisin-2002}.

\begin{theorem}[HLT] The linear map
\[
\aligned
H^{p,q}(X,\bc)&\rightarrow H^{n-q,n-p}(X,\bc)\\
\{\alpha\}&\mapsto \{\alpha\}\smile \{\Omega\}
\endaligned
\]
is an isomorphism, where $\smile$ denotes the cup-product on the cohomology ring $\oplus H^*(X,\bc)$.\end{theorem}

Define the primitive subspace $P^{p,q}(X,\bc)$ of $H^{p,q}(X,\bc)$ by
\[
P^{p,q}(X,\bc):=\Big\{\{\alpha\}\in H^{p,q}(X,\bc),\{\alpha\}\smile\{\Omega\}\smile\{\omega\}=0\Big\}.
\]
Define $Q=Q_\Omega$ a Hermitian form on $H^{p,q}(X,\bc)$ by
\[
Q(\{\alpha\},\{\beta\}):=(\sqrt{-1})^{p-q}(-1)^{\frac{(p+q)(p+q-1)}{2}}\int_X \alpha\wedge\bar{\beta}\wedge\Omega
\] 
for smooth closed $(p,q)$-forms $\alpha$ and $\beta$ on $X$. The integral depends only on the cohomology classes $\{\alpha\}$, $\{\beta\}$ of $\alpha$, $\beta$ in $H^{p,q}(X,\bc)$.

\begin{theorem}[HRR: Hodge-Riemann biliener relations] The Hermitian form $Q$ is positive-definite on $P^{p,q}(X,\bc)$.
\end{theorem}

\begin{theorem}[LD: Lefschetz decomposition theorem] The decomposition
\[
H^{p,q}(X,\bc)=\Big(\{\omega\}\smile H^{p-1,q-1}(X,\bc)\Big)\oplus P^{p,q}(X,\bc)
\]
is orthogonal with respect to the Hermitian form $Q$.
\end{theorem}

Thus we get the signature of $Q$ in terms of the Hodge numbers $h^{p,q}:=\dim H^{p,q}(X,\bc)$. For example when $p=q=1$ we obtain

\begin{Cor}[Hodge index theorem] The signature of $Q$ on $H^{1,1}(X,\bc)$ is $(h^{1,1}-1,1)$.
\end{Cor}

The above theorems are not true if we replace the class of $\Omega=\omega^k$ with an arbitrary class in $H^{k,k}(X,\br)$, even when the class contains a strictly positive form. See Timorin \cite[Remark 3]{Timorin-1998} or Berndtsson-Sibony \cite[Sect. 9]{Berndtsson-Sibony-2002} for an example when $k=2$ and $n=4$.

The sufficient conditions on $\{\Omega\}$ for these theorems have been studied by different groups of mathematicians. We recall some positive results. Let $\omega_1,\dots,\omega_{k+1}$ be arbitrary K\"ahler forms on $X$. Let $\Omega=\omega_1\wedge\omega_2\dots\wedge\omega_k$. Define the primitive space $P^{p,q}(X,\bc)$ and the Hermitian form $Q=Q_\Omega$ with respect to this $\Omega$. More precisely
\[
P^{p,q}(X,\bc):=\Big\{\{\alpha\}\in H^{p,q}(X,\bc),\{\alpha\}\smile\{\Omega\}\smile\{\omega_{k+1}\}=0\Big\}.
\]
and
\[
Q(\{\alpha\},\{\beta\}):=(\sqrt{-1})^{p-q}(-1)^{\frac{(p+q)(p+q-1)}{2}}\int_X \alpha\wedge\bar{\beta}\wedge\Omega.
\] 

\begin{theorem}[mixed HLT, HRR, LD]\label{mixed} For all $p+q\leqslant n$, the class of $\Omega=\omega_1\wedge\dots\wedge\omega_k$ satisfies the hard Lefschetz theorem, the Hodge-Riemann theorem and the Lefschetz decomposition theorem for the bidegree $(p,q)$.
\end{theorem}

Gromov in \cite{Gromov-1990} stated that $Q$ is positive semi-definite when $p=q$. He gave a complete proof for the case $p=q=1$. Later, Timorin in \cite{Timorin-1998} proved the mixed HRR in the linear case, i.e. when $X$ is a complex torus of dimension $n$, see also \cite{Khovanskii-1988, Teissier-1979}. Let $V$ be a complex vector space of dimension $n$ and $\overline{V}$ be its complex conjugate. Let $V^{p,q}=\Lambda^p V\otimes\Lambda^q\overline{V}$ with the convention that $V^{p,q}=0$ unless $0\leqslant p,q\leqslant n$. A form $\omega\in V^{1,1}$ is a {\em K\"ahler form} if $\omega=\sum\limits_{\ell=1}^n\tfrac{\sqrt{-1}}{2}dz_\ell\wedge d\overline{z_\ell}$ in some complex coordinates $(z_1,\dots,z_n)$ of $V$, where $z_i\otimes\overline{z_j}$ is identified with $dz_i\wedge d\overline{z_j}$. Timorin proved the following.

\begin{theorem}[Linear mixed HLT, HRR, LD] Let $\omega_1,\dots,\omega_{k+1}\in V^{1,1}$ be K\"ahler forms. Let $\Omega=\omega_1\wedge\dots\wedge\omega_k$. Then $\wedge\Omega:V^{p,q}\rightarrow V^{n-q,n-p}$ is an isomorphism. Define a Hermitian form on $V^{p,q}$
\[
Q(\alpha,\beta)=(\sqrt{-1})^{p-q}(-1)^{\frac{(p+q)(p+q-1)}{2}}\star(\alpha\wedge\overline{\beta}\wedge\Omega),
\]
where $\star$ is the Hodge-star operator, and define the {\em mixed primitive space}
\[
P^{p,q}=\{\alpha\in V^{p,q}~|~\alpha\wedge\Omega\wedge\omega_{k+1}=0\}
\]
Then $Q$ is positive definite on $P^{p,q}$. Moreover, the decomposition
\[
V^{p,q}=\Big(\omega_{k+1}\wedge V^{p-1,q-1}\Big)\oplus P^{p,q}
\]
is orthogonal with respect to $Q$.
\end{theorem}

Dinh-Nguy\^en in \cite{Dinh-Nguyen-2006} proved Theorem~\ref{mixed} for general compact K\"ahler manifolds, see also Cattani \cite{Cattani-2008} for a later proof using the theory of variations of Hodge structures. Later, in \cite{Dinh-Nguyen-2013} Dinh-Nguy\^en posed a point-wise condition and proved a stronger version.

We will introduce in the next section the notion of Hodge-Riemann cone in the exterior product $V^{k,k}:=\Lambda^k V\otimes\Lambda^k\overline{V}$ with $0\leqslant k\leqslant n$. In practice, $V$ is the complex cotangent space at an arbitrary point $x$ of $X$ and we define Hodge-Riemann cone point-wisely on $X$.

\begin{theorem}[Dinh-Nguy\^en 2013]\label{Dinh-Nguyen-2013} Let $\Omega$ be a closed smooth form of bidegree $(k,k)$ on $X$. Assume that $\Omega$ takes values only in the Hodge-Riemann cone associated with $X$ point-wisely. Then $\{\Omega\}$ satisfies HLT, HRR and LD for all bidegrees $(p,q)$.
\end{theorem}

Roughly speaking, taking values in the Hodge-Riemann cone means $\Omega$ can be continuously deformed to $\omega^k$ in a nice way that some hard Lefschetz properties are preserved. Such deformation does not need to depend on $x$ continuously. Moreover, it does not need to preserve the closedness nor smoothness of the form. The key is to check hard Lefschetz properties point-wisely.

Dinh-Nguy\^en asked in \cite{Dinh-Nguyen-2013} whether the Griffiths cone is contained in the Hodge-Riemann cone. Let $M=(\alpha_{i,j})\in M_{k,k}(V^{1,1})$ be a $k\times k$ matrix with constant coefficient $(1,1)$-form entries. The matrix is {\em Griffiths positive} if $\theta\cdot M\cdot \overline{\theta}^t$ is a K\"ahler form for any $\theta\in\mathbb{C}^k\backslash\{0\}$. The {\em Griffiths cone} is the collection of $(k,k)$-forms which are determinants of Griffiths positive matrices. In particular, for a Griffiths positive vector bundle $(E,h)$ of rank $k$ over a compact complex manifold $X$, the curvature form $\Theta$ is represented by a Griffiths positive $k\times k$ matrix with $(1,1)$-form entries at each point $x\in X$.

Griffiths \cite{Griffiths-1970} proved that the last Chern class of a Griffiths positive vector bundle of rank $2$ is weakly positive. Indeed he proved that the determinant of the curvature matrix is a strictly positive $(2,2)$-form, which implies the positive answer to Dinh-Nguy\^en's question for $n=2$ or $3$.

\begin{Cor}[Answer 1]\label{main1} When $k=2$ and $n=2,3$, the Griffiths cone is contained in the Hodge-Riemann cone.
\end{Cor}

The proof fails when $n\geqslant 4$, since a strictly positive $(2,2)$-form may not satisfy HLT, HRR nor LD. To treat higher dimensional cases, we introduce the following terminology.

Fix complex coordinates $(z_1,\dots,z_n)$ of $\bc^n$. In this paper, to simplify computations, we treat matrices {\em with diagonalized entries}, i.e.
\[
\alpha_{i,j}=\sum\limits_{\ell=1}^n b_{i,j}^{(\ell)}\tfrac{\sqrt{-1}}{2}dz_\ell\wedge d\overline{z_\ell}.
\]
That is to say, there are no terms like $dz_1\wedge d\overline{z_2}$.

\begin{theorem}[Answer 2]\label{main2} Let $M$ be a $k\times k$ matrix with constant coefficient $(1,1)$-form entries. Suppose $M$ is Griffiths positive and has only diagonalized entries. Let $\Omega=\det(M)$. Then
\begin{enumerate}
\item when $n\leqslant 5$, $k=2$, the form $\Omega$ satisfies HLT, HRR and LD for all bidegrees $(p,q)$ with $p+q=n-k$;
\item when $n\geqslant6$, $k=2$, the form $\Omega$ satisfies HLT for bidegrees $(n-2,0)$, $(n-3,1)$, $(1,n-3)$ and $(0,n-2)$.
\end{enumerate}
\end{theorem}

We have two applications dealing forms on compact K\"ahler manifolds, where the condition $M$ having diagonalized entries holds. One referee suggested the first application. When $n\leqslant 5$, it follows from Theorem~\ref{main2}~(1). In fact it is true for all $n$ due to Theorem~\ref{mixed}.

\begin{Cor}\label{app1} Let $X$ be a compact K\"ahler manifold. Let $M$ be a Griffiths positive $2\times 2$ matrix with $(1,1)$-class entries, i.e.
\[
M=
\begin{pmatrix}
\omega_1 & \alpha_{1,2}\\
\overline{\alpha_{1,2}} & \omega_2
\end{pmatrix}.
\]
Suppose that $\alpha_{1,2}=c_1\omega_1+c_2\omega_2$ for some $c_1,c_2\in\mathbb{C}$. Then $\det M$ satisfies HRR, HLT and LD.
\end{Cor}

The second application follows from Theorem~\ref{main2}~(2).

\begin{Cor}\label{app2} Let $X$ be a complex torus of dimension $\leqslant 5$. We identify classes with elements in $V^{p,q}$. Let $M$ be a Griffiths positive $2\times 2$ matrix with diagonalized entries. Then $\det M$ satisfies HRR, HLT and LD.
\end{Cor}

HRR has applications in mixed-volumes inequalities, see Khovanskii \cite{Khovanskii-1978}, Teissier \cite{Teissier-1979,Teissier-1981} and Dinh-Nguy\^en \cite{Dinh-Nguyen-2006}. The reader will find some related results and applications in Cattani \cite{Cattani-2008}, de Cataldo and Migliorini \cite{deCataldo-Migliorini-2002}, Gromov \cite{Gromov-1990}, Dinh and Sibony \cite{Dinh-2012,Dinh-Sibony-2004} and Keum, Oguiso and Zhang \cite{Keum-Oguiso-Zhang-2009,Zhang-2009}. Generalizations to Schur classes are proved by Ross and Toma, see \cite{Ross-Toma-2023-1,Ross-Toma-2023-2,Ross-Toma-2023-3}. Generalization of mixed HRR under $m$-positivity is proved by Xiao in \cite{Xiao-2021}.

This paper is organized as follows. In Section 2 we recall the notion of Hodge-Riemann forms and the Griffiths cone defined by Dinh-Nguy\^en. In Section 3 we explain Corollary~\ref{main1}. In Section 4, we prove Theorem \ref{main2} under the assumption that all entries of $M$ are diagonalized. Section 5 contains two applications.

{\bf Acknowledgement:} The author is grateful to Viet-Anh Nguy\^en (Lille) for introducing the problem, to Weizhe Zheng (MCM, AMSS) for proving Theorem~\ref{hyperdet} about hyperdeterminants and for discussions, and to Baohua Fu (MCM, AMSS) for discussions. The author thanks referees for comments, and for suggesting the applications Corollary~\ref{app1}. The author is supported in part by the Labex CEMPI (ANR-11-LABX-0007-01), the project QuaSiDy (ANR-21-CE40-0016), and China Postdoctoral Science Foundation (2023M733690).

\section{Hodge-Riemann forms}
In this section we recall the Hodge-Riemann form in the linear setting, defined by Dinh-Nguy\^en in \cite{Dinh-Nguyen-2013}.

Let $V$ be an $n$ dimensional complex vector space and $\overline{V}$ be its conjugate space. Denote by
\[
V^{p,q}=\Lambda^p V\otimes\Lambda^q\overline{V}
\]
the space of constant coefficient $(p,q)$-forms with the convention that $V^{p,q}=0$ unless $0\leqslant p,q\leqslant n$. It is a complex vector space of dimension $\binom{n}{p}\binom{n}{q}$. A form $\omega\in V^{1,1}$ is a {\em K\"ahler form} if
\[
\omega=\tfrac{\sqrt{-1}}{2}dz_1\wedge d\overline{z_1}+\dots+\tfrac{\sqrt{-1}}{2}dz_n\wedge d\overline{z_n}
\]
in certain complex coordinates $(z_1,\dots,z_n)$, where $z_i\otimes \overline{z_j}$ is identified with $dz_i\wedge d\overline{z_j}$.

A form $\Omega\in V^{k,k}$ with $0\leqslant k\leqslant n$ is {\em real} if $\Omega=\overline{\Omega}$. Let $V_{\br}^{k,k}$ be the space of all real $(k,k)$-forms. A form $\Omega$ is {\em strictly positive} if its restriction on any $k$ dimentional complex subspace is positive volume form. Fix a K\"ahler form $\omega$ as above.

\begin{definition}[Lefschetz forms]
A $(k,k)$-form $\Omega\in V^{k,k}$ is said to be a {\em Lefschetz form for the bidegree $(p,q)$} if $k=n-p-q$ and the map $\alpha\mapsto \alpha\wedge \Omega$ is an isomorphism between $V^{p,q}$ and $V^{n-q,n-p}$.
\end{definition}

\begin{definition}[Hodge-Riemann forms]
A real $(k,k)$-form $\Omega\in V_{\br}^{k,k}$ is said to be a {\em Hodge-Riemann form for the bidegree $(p,q)$} if there is a continuous deformation $\Omega_t\in V_{\br}^{k,k}$ with $0\leqslant t\leqslant 1$, $\Omega_0=\Omega$ and $\Omega_1=\omega^k$ such that
\[
(*) \ \ \ \ \Omega_t \wedge\omega^{2r}  \ \ \text{is a Lefschetz form for the bidegree $(p-r,q-r)$}
\]
for every $0\leqslant r\leqslant \min\{p,q\}$ and $0\leqslant t\leqslant 1$. The cone of all such forms is called the {\em Hodge-Riemann cone for bidegree $(p,q)$}. We say $\Omega$ is {\em Hodge-Riemann} if it is a Hodge-Riemann form for any bidegree $(p,q)$ with $p+q=n-k$.
\end{definition}

Note that a priori the definition of Hodge-Riemann forms depends on the choice of $\omega$. Due to Dinh-Nguy\^en's results in \cite{Dinh-Nguyen-2013}, a Hodge-Riemann $(k,k)$-form satisfies HLT, HRR and Lefschetz decomposition for all bidegree $(p,q)$ such that $p+q=n-k$. It is proper to use the name ``Hodge-Riemann" here.

According to the classical HLT in the linear case, $(*)$ holds for $t=1$. Moreover, due to Timorin's results in \cite{Timorin-1998}, for any K\"ahler forms $\omega_1,\dots,\omega_k$, the product $\Omega:=\omega_1\wedge\dots\wedge\omega_k$ is a Hodge-Riemann form. In this paper we assume $2\leqslant k\leqslant n$ to avoid the trivial case $k=1$.

Let $M=(\alpha_{i,j})$ be a $k\times k$ matrix with entries in $V^{1,1}$. Assume that $M$ is Hermitian, i.e. $\alpha_{i,j}=\overline{\alpha_{j,i}}$ for all $i,j$. We say that $M$ is {\em Griffiths positive} if for all $\theta=(\theta_1,\dots,\theta_k)\in\bc^k\backslash\{0\}$, $\theta\cdot M\cdot\overline{\theta}^t$ is a K\"ahler form. We call {\em Griffiths cone} the set of $(k,k)$-forms $\Omega:=\det(M)$ with $M$ Griffiths positive.

Dinh-Nguy\^en asked in \cite{Dinh-Nguyen-2013}
\begin{Quest}\label{c2.3} Let $M$ be a Griffiths positive $k\times k$ matrix. Is $\det(M)$ a Hodge-Riemann form?
\end{Quest}

They also explained the relation between this question and
\begin{Quest}\label{c2.4} Let $M$ be a Griffiths positive $k\times k$ matrix. Is $\det(M)$ a Lefschetz form for all bidegree $(p,q)$ such that $p+q=n-k$?
\end{Quest}
in the following sense:

\begin{center}
Yes to Question \ref{c2.3} for $k$ and for some $\omega$ $\Rightarrow$ Yes to Question \ref{c2.4} for $k$.

Yes to Question \ref{c2.4} for $k+2r$, $\forall 0\leqslant r\leqslant (n-k)/2$\\
 $\Rightarrow$ Yes to Question \ref{c2.3} for $k$ and for any $\omega$.
\end{center}

\proof (relation between the two questions) The first implication is due to Dinh-Nguy\^en's Theorem~\ref{Dinh-Nguyen-2013}. Suppose the answer to Question \ref{c2.4} is true for all $k+2r$, $\forall 0\leqslant r\leqslant (n-k)/2$. Then for any Griffiths positive $k\times k$ matrix and for any $t\in[0,1]$, $(1-t)M+t\omega Id_k$ is again a Griffiths positive $k\times k$ matrix. Let $\Omega_t:=\det\big((1-t)M+t\omega {\text Id}_k\big)$. Then $\Omega_t$ is a continuous family such that $\Omega_0=\Omega=\det(M)$ and $\Omega_1=\omega^k$. Moreover, the block matrix
\[
\begin{pmatrix}
(1-t)M+t\omega {\text Id}_k & 0\\
0 & \omega {\text Id}_{2r}
\end{pmatrix}
\]
is a Griffiths positive $(k+2r)\times (k+2r)$ matrix. So its determinant, $\Omega_t\wedge\omega^{2r}$ is Lefschetz for all suitable $p,q,r$. By definition, $\Omega=\det(M)$ is Hodge-Riemann. \qed

In this paper we treat Question \ref{c2.4}, whose statement does not depend on the choice of $\omega$.

\section{Proof for $k=2$ and $n=2,3$}
We cite the appendix of Griffiths \cite{Griffiths-1970}.
\begin{theorem}(Griffiths 1970) Let $E\rightarrow V$ be a Griffiths positive vector bundle with fibre $\bc^2$. Then $c_2(E)>0$.
\end{theorem}

In fact Griffiths has proved that the determinant $\Omega$ of a Griffiths positive $2\times 2$ matrix is strictly positive.

\proof[Proof of Corollary~\ref{main1}] When $n=2$, $\Omega$ is a positive volumn form, hence Lefschetz for bidegree $(0,0)$ and Hodge-Riemann.

When $n=3$, $\Omega$ is a positive volumn form on any $2$ dimensional complex subspace, hence Lefschetz for bidegree $(0,1)$ and $(1,0)$, and Hodge-Riemann. Thus we prove Corollary~\ref{main1}.

This argument does not apply for $k=2,n\geqslant 4$, since Timorin \cite[Remark 3]{Timorin-1998} and  Berndtsson-Sibony \cite[Sect. 9]{Berndtsson-Sibony-2002} already provided strictly positive examples which do not satisfy HLT, HRR nor LD.

\section{Case all entries of $M$ being diagonalized}
In this section, for simplifications, we write $V_j=\tfrac{\sqrt{-1}}{2}dz_j\wedge d\overline{z_j}$ for the euclidean volumn form on the complex line of $z_j$, which is two times the euclidean volumn form. We write $V_{j_1,j_2,\dots,j_s}=V_{j_1}\wedge\cdots\wedge V_{j_s}$ for the volumn form on the $s$ dimensional subspace spanned by those complex lines. Without extra specifications, $\omega=\sum\limits_{\ell=1}^n V_\ell$ denotes the standard K\"ahler form in the linear case. We write $\text{Vol}=V_{1,2,\dots,n}$ for the volumn form on $\bc^n$.

Now we assume that all entries of $M=(\alpha_{i,j})$ are diagonalized, i.e. $\alpha_{i,j}=\sum\limits_{\ell=1}^n b_{i,j}^{(\ell)}V_\ell$. The following lemma holds for general $k\geqslant 2$.

\begin{Prop} Let $M$ be a $k\times k$ matrix with diagonalized entries. We can write $M$ as a matrix valued $(1,1)$-form
\[
M=(\alpha_{i,j})=\Big(\sum\limits_{\ell=1}^n b_{i,j}^{(\ell)}V_\ell\Big)=\sum\limits_{\ell=1}^n (b_{i,j}^{(\ell)})V_\ell
\]
Then $M$ is Griffiths positive if and only if the matrix $B^{(\ell)}:=(b_{i,j}^{(\ell)})$ is a positive definite $k\times k$ matrix for $1\leqslant\ell\leqslant n$.
\end{Prop}
\proof For any $\theta\in\bc^k$, $\theta\neq0$. The $(1,1)$-form
\[
\theta\cdot M\cdot\overline{\theta}^t=\sum\limits_{\ell=1}^n\theta\cdot B^{(\ell)}\cdot\overline{\theta}^tV_\ell.
\]
Thus $M$ is Griffiths positive if and only if $\theta\cdot M\cdot\overline{\theta}^t$ is a K\"ahler form if and only if $\theta\cdot B^{(\ell)}\cdot\overline{\theta}^t>0$ for $1\leqslant\ell\leqslant n$.\qed

After a dilation on coordinates we may assume that $\alpha_{1,1}=\omega=\sum\limits_{\ell=1}^n V_\ell$, i.e. $b_{1,1}^{(\ell)}=1$ for $1\leqslant \ell\leqslant n$.

To calculate $\Omega=\det(M)$, we introduce the {\em hyperdeterminant} among $B^{(\ell)}$.

\begin{definition} Let $B^{(1)},\dots,B^{(k)}$ be $k\times k$ complex valued matrices with $B^{(\ell)}=(b^{(\ell)}_{i,j})$. We define the $k\times k\times k$ hypermatrix ${\bf B}=(B^{(1)},\dots,B^{(k)})$ a $3$ dimensional array whose layers $B^{(\ell)}$ are matrices. We define its {\em hyperdeterminant} by
\[
{\sf hdet}({\bf B}):=\sum\limits_{\sigma,\tau\in S_k}\text{sgn}(\sigma)\prod\limits_{j=1}^k b_{j,\sigma(j)}^{(\tau(j))}
\]
where $S_k$ is the permutation group of $k$ elements.
\end{definition}

We remark that switching two layers does not change the hyperdeterminant. The determinant
\[
\Omega=\det(M)=\sum\limits_{1\leqslant i_1<\dots<i_k\leqslant n}\Omega_{i_1,\dots,i_k} V_{i_1,\dots,i_k}
\]
where $\Omega_{i_1,\dots,i_k}={\sf hdet}\big((B^{(i_1)},\dots,B^{(i_k)})\big)$, whose positivity is proved by Weizhe Zheng (MCM, AMSS).
\begin{theorem}[Zheng 2021]\label{hyperdet} Let $B^{(1)},\dots,B^{(k)}$ be positive semidefinite Hermitian matrices. Denote by $\mu^{(\ell)}$ and $\lambda^{(\ell)}$ the minimal and the maximal eigenvalues of $B^{(\ell)}$. Then
\[
k!\mu^{(1)}\cdots\mu^{(k)}\leqslant {\sf hdet}({\bf B})\leqslant k!\lambda^{(1)}\cdots\lambda^{(k)}.
\]
In particular, if $B^{(1)},\dots,B^{(k)}$ are positive definite, then ${\sf hdet}({\bf B})>0$.
\end{theorem}

\proof We proceed by induction on $k$. The case $k=1$ is trivial.

For the general case, note that for any $k\times k$ matrix $U$, we have ${\sf hdet}({\bf B}U)=\det(U){\sf hdet}({\bf B})={\sf hdet}(U{\bf B})$, where ${\bf B}U$ and $U{\bf B}$ are the hypermatrices defined by layerwise multiplication $({\bf B}U)^{(\ell)}:=B^{(\ell)}U$ and $(U{\bf B})^{(\ell)}:=UB^{(\ell)}$. Thus, up to replacing ${\bf B}$ by $U{\bf B}U^H$ for a unitary matrix $U$, we may assume that one layer, say $B^{(1)}$, is diagonal.

Once $B^{(1)}$ is diagonal, in the definition of ${\sf hdet}({\bf B})$ above, if $\tau(j)=1$, then $b_{j,\sigma(j)}^{(1)}$ is nonzero only if $\sigma(j)=j$. Thus
\[
{\sf hdet}({\bf B})=\sum\limits_{j=1}^k b^{(1)}_{j,j}{\sf hdet}({\bf B}_j)
\]
where ${\bf B}_j$ is the $(k-1)\times(k-1)\times(k-1)$ hypermatrix obtained from ${\bf B}$ by removing the layer $B^{(1)}$ and removing the $j$-th row together with the $j$-th column in each layer $B^{(i)}$. We conclude by the induction hypothesis. \qed

\begin{Cor} Let $M$ be a $k\times k$ Griffiths positive matrix with diagonalized entries. Then $\Omega=\det(M)$ is a strictly positive $(k,k)$-form, hence a Lefschetz form for bidegree $(n-k,0)$ and $(0,n-k)$.
\end{Cor}

\begin{Cor} Let $M$ be a $k\times k$ Griffiths positive matrix with diagonalized entries. If $k=n-1$ or $k=n$, then $\det(M)$ satisfies HLT, HRR, LD for all suitable bidegrees.
\end{Cor}

\subsection{Case $k=2$, $n=4$, $M$ with diagonalized entries}
When $k=2$, each positive definite matrix $B^{(\ell)}$ can be written as
\[
B^{(\ell)}=
\begin{pmatrix}
1 &  b_\ell\\
\overline{b_\ell} & |b_\ell|^2+t_\ell
\end{pmatrix}
\]
for some $b_\ell\in\bc$ and $t_\ell>0$. The determinant
\[
\Omega=\det(M)=\sum\limits_{1\leqslant i<j\leqslant 4}\Omega_{i,j}V_{i,j}
\]
where
\[
\Omega_{i,j}={\sf hdet}(B^{(i)},B^{(j)})=|b_i|^2+t_i+|b_j|^2+t_j-b_i\overline{b_j}-b_j\overline{b_i}=|b_i-b_j|^2+t_i+t_j>0.
\]
For simplifications, we define $b_{i,j}:=b_i-b_j\in\bc$.

\begin{theorem}\label{thm-n4k2} Let $n=4$. Let $M$ be a $2\times 2$ Griffiths positive matrix with diagonalized entries. Then $\Omega=\det(M)$ is a Lefschetz form for bidegree $(2,0)$, $(1,1)$ and $(0,2)$.
\end{theorem}
\proof It suffices to check bidegree $(2,0)$ and $(1,1)$.

(Trivial part) For bidegree $(2,0)$, we take the standard basis of $V^{2,0}$ by the lexicographical order
\[
\{dz_1\wedge dz_2,dz_1\wedge dz_3,\dots,dz_3\wedge dz_4\}
\]
and we take a basis of $V^{4,2}$ by the Hodge-star of their conjugates
\[
\{dz_1\wedge dz_2\wedge V_{3,4},dz_1\wedge dz_3\wedge V_{2,4},\dots,dz_3\wedge dz_4\wedge V_{1,2}\}.
\]
Under these two basis, the linear map $\wedge\Omega:V^{2,0}\rightarrow V^{4,2}$ can be expressed as a diagonal matrix $\text{diag}\{\Omega_{3,4},\Omega_{2,4},\dots,\Omega_{1,2}\}$ with positive entries. This map is an isomorphism.

(Non-trivial part) For bidegree $(1,1)$, again, we take the standard basis of $V^{1,1}$ as follows
\[
\{\underbrace{V_1,V_2,V_3,V_4}_{\text{$4$ elements}},\underbrace{dz_1\wedge d\overline{z_2},dz_1\wedge d\overline{z_3},\dots,dz_4\wedge d\overline{z_3}}_{\text{$12$ elements}}\}
\]
and we take a basis of $V^{3,3}$ by the Hodge star of their conjugates
\[
\{
\underbrace{V_{2,3,4},V_{1,3,4},V_{1,2,4},V_{1,2,3}}_{\text{$4$ elements}},\underbrace{dz_1\wedge d\overline{z_2}\wedge V_{3,4},dz_1\wedge d\overline{z_3}\wedge V_{2,4},\dots,dz_4\wedge d\overline{z_3}\wedge V_{1,2}}_{\text{$12$ elements}}.
\}
\]

Under these two basis, the linear map $\wedge \Omega:V^{1,1}\rightarrow V^{3,3}$ can be expressed as a blocked matrix
\[
\begin{pmatrix}
G &  0\\
0 & \text{diag}\{\Omega_{3,4},\Omega_{2,4},\dots,\Omega_{1,2}\}
\end{pmatrix}
\]
where
\[
G=
\begin{pmatrix}
 0 & \Omega_{3,4} & \Omega_{2,4} & \Omega_{2,3} \\
\Omega_{3,4} &  0 & \Omega_{1,4} & \Omega_{1,3} \\
\Omega_{2,4} & \Omega_{1,4} & 0 & \Omega_{1,2} \\
\Omega_{2,3} & \Omega_{1,3} & \Omega_{1,2} &  0
\end{pmatrix}.
\]
It suffices to verify that $\det(G)\neq0$. Let
\[
A=\sqrt{\Omega_{3,4}\Omega_{1,2}}, \ \ B=\sqrt{\Omega_{2,4}\Omega_{1,3}}, \ \ C=\sqrt{\Omega_{2,3}\Omega_{1,4}}.
\]
Then $\det(G)=-(A+B+C)(A+B-C)(A-B+C)(-A+B+C)$ has the form of the Heron formula which calculates the area of a triangle with side length $(A,B,C)$. We are going to show that the side lengths $(A,B,C)$ actually forms a triangle. After permutations among $b_j$ and among $t_j$, it suffices to verify that
\[
\sqrt{\Omega_{3,4}\Omega_{1,2}}+\sqrt{\Omega_{2,4}\Omega_{1,3}}>\sqrt{\Omega_{2,3}\Omega_{1,4}},
\]
i.e.
\[
\Omega_{3,4}\Omega_{1,2}+\Omega_{2,4}\Omega_{1,3}+2\sqrt{\Omega_{3,4}\Omega_{1,2}}\sqrt{\Omega_{2,4}\Omega_{1,3}}-\Omega_{2,3}\Omega_{1,4}>0
\]
The left hand side is
\[
\aligned
LHS=  & (|b_{3,4}|^2+t_3+t_4)(|b_{1,2}|^2+t_1+t_2)+
(|b_{2,4}|^2+t_2+t_4)(|b_{1,3}|^2+t_1+t_3)\\
 & + 2\sqrt{(|b_{3,4}|^2+t_3+t_4)(|b_{1,2}|^2+t_1+t_2)(|b_{2,4}|^2+t_2+t_4)(|b_{1,3}|^2+t_1+t_3)}\\
& -(|b_{2,3}|^2+t_2+t_3)(|b_{1,4}|^2+t_1+t_4)
\endaligned
\]
We expand the product in the square-root and sort each summand according to the order of $t_j$:
\[
\Omega_{3,4}\Omega_{1,2}\Omega_{2,4}\Omega_{1,3}=\underbrace{(\dots)}_{\text{$b$ part}}+\sum\limits_{j=1}^4\underbrace{(\dots)t_j}_{\text{$t_j$ part}}+\sum\limits_{j=1}^4\underbrace{(\dots)t_j^2}_{\text{$t_j^2$ part}}+\sum\limits_{1\leqslant j<\ell\leqslant 4}^4\underbrace{(\dots)t_jt_\ell}_{\text{$t_j t_\ell$ part}}+O_{t_j}(3)
\]
where
\[
\aligned
\text{($b$ part)}&:=|b_{3,4}|^2|b_{1,2}|^2|b_{2,4}|^2|b_{1,3}|^2\\
\text{($t_1^2$ part)}&:=\big(t_1|b_{3,4}||b_{2,4}|\big)^2\\
\text{($t_1$ part)}&:=t_1|b_{3,4}|^2|b_{2,4}|^2\big(|b_{1,2}|^2+|b_{1,3}|^2\big)\\
&\geqslant 2t_1|b_{3,4}|^2|b_{2,4}|^2|b_{1,2}||b_{1,3}|\\
&=2\sqrt{\text{($t_1^2$ part)}}\sqrt{\text{($b$ part)}}\\
\text{($t_1t_2$ part)}&:=t_1t_2|b_{3,4}|^2\big(|b_{1,2}|^2+|b_{2,4}|^2+|b_{1,3}|^2\big)\\
&\geqslant 2t_1t_2|b_{3,4}|^2|b_{2,4}||b_{1,3}|\\
&=2\sqrt{\text{($t_1^2$ part)}}\sqrt{\text{($t_2^2$ part)}}\\
\text{($t_1t_4$ part)}&:=t_1t_4\big(|b_{2,4}|^2+|b_{3,4}|^2\big)\big(|b_{1,2}|^2+|b_{1,3}|^2\big)\\
&\geqslant 2t_1t_4|b_{1,2}||b_{1,3}||b_{2,4}||b_{3,4}|\\
&=2\sqrt{\text{($t_1^2$ part)}}\sqrt{\text{($t_4^2$ part)}}
\endaligned
\]
In fact $\text{($t_j$ part)}\geqslant 2\sqrt{\text{($t_j^2$ part)}}\sqrt{\text{($b$ part)}}$ and $\text{($t_jt_\ell$ part)}\geqslant 2\sqrt{\text{($t_j^2$ part)}}\sqrt{\text{($t_\ell^2$ part)}}$
for any $j,l$. Hence
\[
\aligned
2&\sqrt{\Omega_{3,4}\Omega_{1,2}\Omega_{2,4}\Omega_{1,3}}\\
>2&\sqrt{\text{($b$ part)}+\sum\limits_{j=1}^4\text{($t_j^2$ part)}+\sum\limits_{j=1}^4 2\sqrt{\text{($t_j^2$ part)}}\sqrt{\text{($b$ part)}}+\sum\limits_{1\leqslant j<l\leqslant 4}2\sqrt{\text{($t_j^2$ part)}}\sqrt{\text{($t_\ell^2$ part)}}}\\
=2&\Big(\sqrt{\text{($b$ part)}}+\sum\limits_{j=1}^4 \sqrt{\text{($t_j^2$ part)}}\Big).
\endaligned
\]
Here the first inequality is strict because there is a term $t_1t_2t_3t_4>0$. Thus
\[
\aligned
LHS>&|b_{3,4}|^2|b_{1,2}|^2+|b_{2,4}|^2|b_{1,3}|^2+2\sqrt{\text{($b$ part)}}-|b_{2,3}|^2|b_{1,4}|^2\\
&+t_1\big(|b_{3,4}|^2+|b_{2,4}|^2+\underbrace{2|b_{3,4}||b_{2,4}|}_{\text{comes from the square root}}-|b_{2,3}|^2\big)\\
&+t_2\big(|b_{3,4}|^2+|b_{1,3}|^2+2|b_{3,4}||b_{1,3}|-|b_{1,4}|^2\big)\\
&+t_3\big(|b_{1,2}|^2+|b_{2,4}|^2+2|b_{1,2}||b_{2,4}|-|b_{1,4}|^2\big)\\
&+t_4\big(|b_{1,2}|^2+|b_{1,3}|^2+2|b_{1,2}||b_{1,3}|-|b_{2,3}|^2\big)\\
&+\underbrace{(t_3+t_4)(t_1+t_2)+(t_2+t_4)(t_1+t_3)-(t_2+t_3)(t_1+t_4)}_{=0}.
\endaligned
\]
The first line is $(|b_{3,4}||b_{1,2}|+|b_{2,4}||b_{1,3}|)^2-|b_{2,3}|^2|b_{1,4}|^2$. In fact
\[
\aligned
b_{2,3}b_{1,4}&=(b_2-b_3)(b_1-b_4)\\
&=(b_3-b_4)(b_2-b_1)+(b_2-b_4)(b_1-b_3)=-b_{3,4}b_{1,2}+b_{2,4}b_{1,3}\\
|b_{2,3}b_{1,4}|&\leqslant |b_{3,4}b_{1,2}|+|b_{2,4}b_{1,3}|
\endaligned
\]
Indeed this is Ptolemy Theorem.

\begin{figure}[htbp]
\begin{center}
   \includegraphics[width=0.3\linewidth]{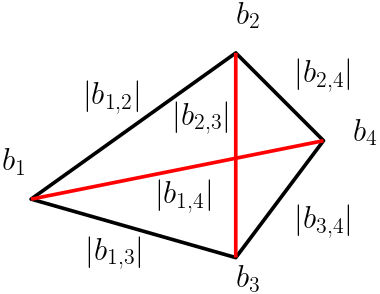}
   \end{center}
\end{figure}

In the second line
\[
t_1\big(|b_{3,4}|^2+|b_{2,4}|^2+2|b_{3,4}||b_{2,4}|-|b_{2,3}|^2\big)=t_1\big((|b_{3,4}|+|b_{2,4}|)^2-|b_{2,3}|^2\big)\geqslant 0
\]
by the triangle inequality. So $LHS>0$ which implies $A+B-C>0$. After permutations among $b_j$ and among $t_j$, we conclude that $\det(G)<0$ and $\wedge \Omega:V^{1,1}\rightarrow V^{3,3}$ is an isomorphism.\qed

\subsection{Case $k=2$, $n\geqslant 4$, $M$ with diagonalized entries}
As before, write $M=\sum\limits_{\ell=1}^n B^{(\ell)} V_\ell$ with
\[
B^{(\ell)}=
\begin{pmatrix}
1 &  b_\ell\\
\overline{b_\ell} & |b_\ell|^2+t_\ell
\end{pmatrix}
\]
for some $b_\ell\in\bc$ and $t_\ell>0$. The determinant
\[
\Omega=\det(M)=\sum\limits_{1\leqslant i<j\leqslant n}\Omega_{i,j}V_{i,j}
\]
where
\[
\Omega_{i,j}=|b_{i,j}|^2+t_i+t_j>0.
\]

\begin{theorem} Let $n\geqslant 4$. Let $M$ be a $2\times 2$ Griffiths positive matrix with diagonalized entries. Then $\Omega=\det(M)$ is a Lefschetz form for bidegree $(n-2,0)$, $(n-3,1)$, $(1,n-3)$ and $(0,n-2)$.
\end{theorem}
\proof The technique is the same as in Theorem \ref{thm-n4k2}. We only need to choose basis carefully.

For bidegree $(n-2,0)$, take the lexicographical ordered basis
\[
\{dz_1\wedge dz_2\wedge\dots\wedge dz_{n-2},\dots,dz_3\wedge\dots\wedge dz_n\}
\]
of $V^{n-2,0}$ and take the Hodge-star of their conjugates as basis of $V^{n,2}$. Then the matrix form of the linear map $\wedge\Omega:V^{n-2,0}\rightarrow V^{n,2}$ is diag$(\Omega_{n-1,n},\Omega_{n-2,n},\dots,\Omega_{1,2})$ where the indices are in the reversed lexicographical order. Each $\Omega_{i,j}>0$ implies that $\Omega$ is a Lefschetz form for bidegree $(n-2,0)$.

For bidegree $(n-3,1)$, take basis of $V^{n-3,1}$ as follows
\[
\aligned
\{
&dz_1\wedge\dots\wedge dz_{n-4}\wedge V_{n-3}, \quad dz_1\wedge\dots\wedge dz_{n-4}\wedge V_{n-2},\\
&dz_1\wedge\dots\wedge dz_{n-4}\wedge V_{n-1}, \quad dz_1\wedge\dots\wedge dz_{n-4}\wedge V_{n},\\
&dz_1\wedge\dots\wedge dz_{n-4}\wedge dz_{n-3}\wedge d\overline{z_{n-2}}, \quad \dots, \quad dz_1\wedge\dots\wedge dz_{n-4}\wedge dz_n\wedge d\overline{z_{n-1}},\\
&dz_1\wedge\dots\wedge dz_{n-5}\wedge dz_{n-3}\wedge V_{n-4}, \quad \dots, \quad dz_1\wedge\dots\wedge dz_{n-5}\wedge dz_{n-3}\wedge V_{n},\\
&dz_1\wedge\dots\wedge dz_{n-5}\wedge dz_{n-3}\wedge (\text{$(1,1)$-forms s.t. the product contains no $V_j$}), \quad \dots,\\
&dz_5\wedge\dots\wedge dz_n\wedge V_1, \quad \dots, \quad 
dz_5\wedge\dots\wedge dz_n\wedge V_4,\\
&dz_5\wedge\dots\wedge dz_n\wedge dz_1\wedge d\overline{z_2}, \quad \dots, \quad dz_5\wedge\dots\wedge dz_n\wedge dz_4\wedge d\overline{z_3}
\}
\endaligned
\]
and take basis of $V^{n-1,3}$ by the Hodge-star of their conjugates. Then the matrix form of the linear map $\wedge\Omega:V^{n-3,1}\rightarrow V{n-1,3}$ is
\[
\begin{pmatrix}
G_{n-3,n-2,n-1,n} &  0 &  \dots & 0 & 0\\
0 & \text{diag}\{\Omega_{n-1,n},\dots,\Omega_{n-3,n-2}\} &  \dots & 0 & 0\\
\vdots & \vdots &  \ddots & \vdots & \vdots \\
0 & 0   & \dots & G_{1,2,3,4} & 0\\
0 & 0  & \dots & 0 & \text{diag}\{\Omega_{3,4},\dots,\Omega_{1,2}\}
\end{pmatrix}
\]
where each
\[
G_{i_1,i_2,i_3,i_4}=
\begin{pmatrix}
 0 & \Omega_{i_3,i_4} & \Omega_{i_2,i_4} & \Omega_{i_2,i_3} \\
\Omega_{i_3,i_4} &  0 & \Omega_{i_1,i_4} & \Omega_{i_1,i_3} \\
\Omega_{i_2,i_4} & \Omega_{i_1,i_4} & 0 & \Omega_{i_1,i_2} \\
\Omega_{i_2,i_3} & \Omega_{i_1,i_3} & \Omega_{i_1,i_2} &  0
\end{pmatrix}
\]
is invertible. So $\Omega$ is a Lefschetz form for bidegree $(n-3,1)$.\qed

\begin{Cor} Let $n=4,5$. Let $M$ be a $2\times 2$ Griffiths positive matrix with diagonalized entries. Then $\det(M)$ is a Hodge-Riemann form.
\end{Cor}

\subsection{Difficulty of the case $k=2$, $n=6$ and $M$ with diagonalized entries}
We already checked that $\Omega$ is a Lefschetz form for bidegree $(4,0)$, $(3,1)$, $(1,3)$ and $(0,4)$. The only thing left is the bidegree $(2,2)$. It amounts to prove that the matrix

\begin{figure}[htbp]
\begin{center}
   \includegraphics[width=\linewidth]{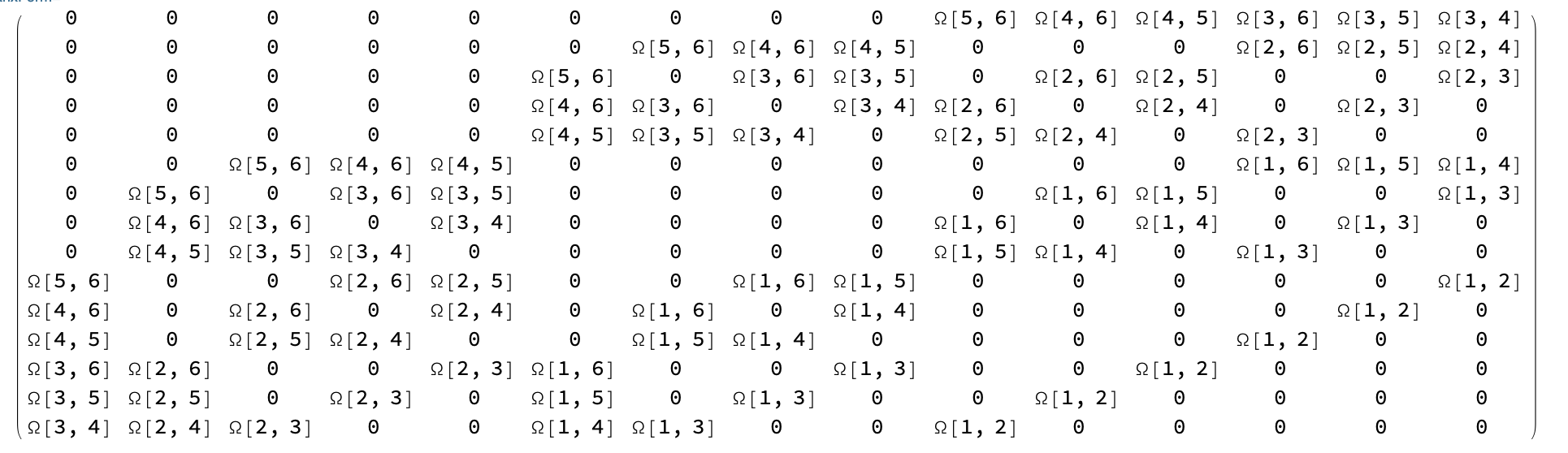}
   \end{center}
\end{figure}

\noindent is invertible. However, unlike the case before, the determinant is irreducible in $\bc[\sqrt{\Omega_{i,j}}]$. Indeed, by a Mathematica program, the determinant is irreducible on $\bc[\Omega_{i,j}^{1/t}]$ for $t=1,2,\dots,15$. It is difficult to prove that this determinant is non-zero by the techniques, analogues of Heron formula, as before.

\section{Applications on compact K\"ahler manifolds and on torus}
Let $X$ be a compact K\"ahler manifold dimension $n$. Let $M$ be a Griffiths positive $2\times 2$ matrix with entires in $H^{1,1}(X)$
\[
M=
\begin{pmatrix}
\alpha_{1,1} & \alpha_{1,2}\\
\alpha_{2,1} & \alpha_{2,2}
\end{pmatrix}.
\]
The Griffiths positivity implies that $\alpha_{1,1}=:\omega_1$ and $\alpha_{2,2}=:\omega_2$ are K\"ahler, and $\alpha_{1,2}=\overline{\alpha_{2,1}}$. So we may write
\[
M=
\begin{pmatrix}
\omega_1 & \alpha_{1,2}\\
\overline{\alpha_{1,2}} & \omega_2
\end{pmatrix}.     
\]

In this section we put the following assumption: {\em all entries of $M$ are linear combinations of the diagonal entries}, i.e. $\alpha_{1,2}=c_1\omega_1+c_2\omega_2$ for some $c_1,c_2\in\mathbb{C}$.

\begin{lemma} If $\alpha_{1,2}=c_1\omega_1+c_2\omega_2$, then at each point of $X$, there are holomorphic coordinates such that $M$ has only diagnolaized entires.
\end{lemma}
\proof At each point of $X$, there are holomorphic coordinates $(z_1,\dots,z_n)$ such that $\omega_1=\sum\limits_{i=1}^n \frac{\sqrt{-1}}{2}d z_i\wedge d\overline{z_i}$. After a unitary change of coordiantes we may assume $\omega_2=\sum\limits_{i=1}^n\lambda_i\frac{\sqrt{-1}}{2}d z_i\wedge d\overline{z_i}$ for some $\lambda_i>0$, while the expression of $\omega_1$ is invariant. Then $\alpha_{1,2}=\sum\limits_{i=1}^n(c_1+c_2\lambda_i)\frac{\sqrt{-1}}{2}d z_i\wedge d\overline{z_i}$ and $\overline{\alpha_{1,2}}$ are diagonalized. \qed

In this section we prove Theorem~\ref{app1} that $\det M$ satisfies HRR, HLT, LD. When $n\leqslant 5$, it follows directly from Theorem~\ref{main2}. But indeed it is true in any dimension. The idea is that $\det M=\tilde{\omega_1}\wedge\tilde{\omega_2}$ can be written as a product of two new K\"ahler classes and the result follows from Theorem~\ref{mixed}.

\begin{lemma}\label{lem5.2} If $\alpha_{1,2}=c_1\omega_1+c_2\omega_2$, then $M$ is Griffiths positive if and only if 
\[
|c_1+\lambda_{min}c_2|^2<\lambda_{min}, \ \ \ \ |c_1+\lambda_{max}c_2|^2<\lambda_{max},
\]
where
\[
\lambda_{min}=\sup\{\lambda>0~|~\omega_2-\lambda\omega_1~\text{is K\"ahler}\}, \ \ \ \ 
\lambda_{max}=\inf\{\lambda>0~|~\lambda\omega_1-\omega_2~\text{is K\"ahler}\}.
\]
\end{lemma}

\proof At each point $z\in M$, under certain coordinates, the entires can be expressed as matrices
\[
\omega_1=I_n,  \  \  \  \  \omega_2=\text{diag}\{\lambda_1(z),\dots,\lambda_n(z)\}
\]
for some $\lambda_1(z),\dots,\lambda_n(z)>0$. Thus
\[
\lambda_{min}=\min\{\lambda_i(z)~|~z\in M,1\leqslant i\leqslant n\}, \ \ \ \ 
\lambda_{max}=\max\{\lambda_i(z)~|~z\in M,1\leqslant i\leqslant n\}.
\]

$M$ is Griffiths positive if and only if for each $\zeta\in\mathbb{C}$
\[
(1,\zeta)\cdot M\cdot(1,\overline{\zeta})^t=\text{diag}\{1+c_1\zeta+\lambda_i(z)c_2\zeta+\overline{c_1\zeta}+\lambda_i(z)\overline{c_2\zeta}+\lambda_i(z)|\zeta|^2\}_{i=1}^n
\]
is K\"ahler, i.e. $\forall z\in M$, $ 1\leqslant i\leqslant n$, $\zeta\in\mathbb{C}$,
\begin{align*}
&
1+c_1\zeta+\lambda_i(z)c_2\zeta+\overline{c_1\zeta}+\lambda_i(z)\overline{c_2\zeta}+\lambda_i(z)|\zeta|^2
\\
=&\left|\sqrt{\lambda_i(z)}\zeta+\frac{\overline{c_1+\lambda_i(z)c_2}}{\sqrt{\lambda_i(z)}}\right|^2+1-\frac{|c_1+\lambda_i(z)c_2|^2}{\lambda_i(z)}>0,
\end{align*}
i.e.
\[
|c_1+\lambda_i(z)c_2|^2<\lambda_i(z).
\]

Fix $c_1,c_2\in\mathbb{C}$, the function $f:\mathbb{R}\rightarrow\mathbb{R}$, $f(t):=|c_1+t\,c_2|^2-t$ is a convex, i.e. $f''(t)\geqslant 0$. Hence 
$f\big(\lambda_i(z)\big)<0$ for each $1\leqslant i\leqslant n$ and for all $z\in M$ if and only if $f(\lambda_{min})<0$ and $f(\lambda_{max})<0$. \qed

\proof[Proof of Corollary~\ref{app1}] The case $\lambda_{min}=\lambda_{max}$ is trivial, since all entries are proportional, $\det(M)$ is proportional to $\omega_1\wedge\omega_1$. The case $c_1\,c_2=0$ is also trivial, since
\[
\det
\begin{pmatrix}
\omega_1 & c_1\,\omega_1\\
\overline{c_1}\,\omega_1 & \omega_2
\end{pmatrix}
=
\omega_1\wedge(\omega_2-|c_1|^2\,\omega_1)
\]
is a product of two K\"ahler classes.

Now we assume $\lambda_{min}<\lambda_{max}$ and $c_1\,c_2\neq 0$. We want to write $\det(M)=\tilde{\omega_1}\wedge\tilde{\omega_2}$, where
\[
\tilde{\omega_1}=-\frac{u}{|c_2|^2}\,\omega_1+\omega_2, \ \ \ \ \tilde{\omega_2}=v\,\omega_1-|c_2|^2\omega_2,
\]
for some $u,v\in\mathbb{R}$ to be determined. The forms $\tilde{\omega_1}$ and $\tilde{\omega_2}$ are K\"ahler if and only if
\begin{align}\label{interval}
\left\{
\begin{aligned}
0<u&<\lambda_{min}|c_2|^2\\
v&>\lambda_{max}|c_2|^2
\end{aligned}
\right.
\end{align}
Expanding
\[
\det(M)=-|c_1|^2\omega_1^2-|c_2|^2\omega_2^2+(1-c_1\overline{c_2}-\overline{c_1}c_2)\omega_1\wedge\omega_2,
\]
then
\[
\left\{
\begin{aligned}
u\,v&=|c_1\,c_2|^2\\
u+v&=1-c_1\overline{c_2}-\overline{c_1}c_2
\end{aligned}
\right.
\]
The system is solvable if and only if $u,v$ are exactly the two roots of the quadratic function $f(x):=x^2-(1-c_1\overline{c_2}-\overline{c_1}c_2)+|c_1\,c_2|^2$ lies in the intervals (\ref{interval}), if and only if
\[
\left\{
\begin{aligned}
f(\lambda_{min}|c_2|^2)=\lambda_{min}^2|c_2|^4-(1-c_1\overline{c_2}-\overline{c_1}c_2)\lambda_{min}|c_2|^2+|c_1\,c_2|^2&<0\\
f(\lambda_{max}|c_2|^2)=\lambda_{max}^2|c_2|^4-(1-c_1\overline{c_2}-\overline{c_1}c_2)\lambda_{max}|c_2|^2+|c_1\,c_2|^2&<0
\end{aligned}
\right.
\]

By Lemma~\ref{lem5.2}, the Griffiths positivity condition is equivalent to
\[
|c_1+\lambda_{min}c_2|^2-\lambda_{min}<0, \ \ \ \ |c_1+\lambda_{max}c_2|^2-\lambda_{max}<0.
\]
Expanding the sqaures
\[
\aligned
\lambda_{min}^2|c_2|^2-(1-c_1\overline{c_2}-\overline{c_1}c_2)\,\lambda_{min}+|c_1|^2<0,\\
\lambda_{max}^2|c_2|^2-(1-c_1\overline{c_2}-\overline{c_1}c_2)\,\lambda_{max}+|c_1|^2<0.
\endaligned
\]
Thus
\[
\aligned
f(\lambda_{min}|c_2|^2)<0, \ \ \ \ f(\lambda_{max}|c_2|^2)<0.
\endaligned
\]
implies that $\det(M)$ is a product of two K\"ahler classes. The proof finished by applying Theorem~\ref{mixed}.\qed

\proof[Proof of Corollary~\ref{app2}]
Let $n\leqslant 5$, $\Lambda\subset\mathbb{C}^n$ be a lattice of rank $2n$ and let $X=\mathbb{C}^n/\Lambda$ be a complex torus of dimension $n$. There is a nature ring isomorphism between the cohomology ring $\bigoplus_{p,q}H^{p,q}(X,\mathbb{C})$ and $\bigoplus_{p,q}\Lambda^p(\mathbb{C}^n)\otimes\Lambda^p(\mathbb{C}^n)$. Let $M$ be the matrix with $(1,1)$-class entries
\[
M=
\begin{pmatrix}
\omega_1 & \alpha_{1,2}\\
\overline{\alpha_{1,2}} & \omega_2
\end{pmatrix}
\]
where entries can be simultaneously diagnolized. After a linear change of coordinate we may write
\[
\omega_1=I_n,  \ \ \ \ \alpha_{1,2}=\text{diag}\{b_i\}_{i=1}^n, \ \ \ \ \omega_2=\text{diag}\{|b_i|^2+t_i\}_{i=1}^n,
\]
for some $b_i\in\mathbb{C}$ and $t_i>0$. According to Theorem~\ref{main2}, the $(2,2)$-class
\[
\det(M)=\sum\limits_{1\leqslant i<j\leqslant n}\Omega_{i,j}\,V_{i,j}, \ \ \ \ \text{where~}\Omega_{i,j}=|b_i-b_j|^2+t_i+t_j,
\]
satisfies HRR, HLT and LD.\qed

It is still not clear whether this kind of $(2,2)$-class can be written as a wedge product of two K\"ahler classes. In general $\det(M)$ cannot be written as a product of two diagonalized K\"ahler classes
\[
\tilde{\omega}_1=\text{diag}\{\lambda_{1,i}\}_{i=1}^n, \ \ \ \ \tilde{\omega}_2=\text{diag}\{\lambda_{2,i}\}_{i=1}^n;
\]
for some $\lambda_{1,i},\lambda_{2,i}>0$, since
\[
\tilde{\omega}_1\wedge\tilde{\omega}_2=\sum\limits_{1\leqslant i<j\leqslant n}\tilde{\Omega}_{i,j}\,V_{i,j}
\]
must satisfy $\tilde{\Omega}_{i,j}\,\tilde{\Omega}_{i',j'}=\tilde{\Omega}_{i,j'}\,\tilde{\Omega}_{i',j}$. But the coefficients $\Omega_{i,j}$ in $\det(M)$ do not satisfy
$\Omega_{1,2}\,\Omega_{3,4}=\Omega_{1,4}\,\Omega_{2,3}$ in general, e.g. when $b_j=0$ and $t_j=j$ for $j\in\{0,\dots,n\}$.

\begin{center}
	\bibliographystyle{amsplain}
	\bibliography{references}
\end{center}
\setlength\parindent{0em}
{\scriptsize Zhangchi Chen, Morningside Center of Mathematics, Academy of Mathematics and Systems Science, Chinese Academy of Sciences, Beijing 100190, China}\\
{\bf\scriptsize zhangchi.chen@amss.ac.cn}, {\bf\scriptsize http://www.mcm.ac.cn/people/postdocs/202110/t20211022\_666685.html}

\end{document}